\documentclass[12pt]{amsart}

\usepackage{amssymb}

\usepackage{enumitem}

\usepackage{graphicx}

\makeatletter
\@namedef{subjclassname@2020}{%
  \textup{2020} Mathematics Subject Classification}
\makeatother

\usepackage[T1]{fontenc}

\newtheorem{theorem}{Theorem}[section]
\newtheorem{corollary}[theorem]{Corollary}
\newtheorem{lemma}[theorem]{Lemma}
\newtheorem{proposition}[theorem]{Proposition}

\newtheorem{question}[theorem]{Question}

\newtheorem{claim}{Claim}

\theoremstyle{definition}

\newtheorem{example}[theorem]{Example}

\numberwithin{equation}{section}

\usepackage[hidelinks]{hyperref}

\newcommand{\abs}[1]{\lvert#1\rvert}
\newcommand{\res}{\upharpoonright}

\DeclareMathOperator{\Clop}{Clop}

\begin{document}

\title{Rigidity of ordinal spaces}


\author[Yi Zhou]{Yi Zhou}
\address{School of Sciences, Southwest Petroleum University\\
Chengdu, 610500, China}
\email{yizhou@swpu.edu.cn}

\author[Jiaheng Zhuang]{Jiaheng Zhuang$^{*}$}
\address{School of Mathematics, Sichuan University\\
Chengdu, 610064, China}
\email{jiahengzhuang@outlook.com}
\thanks{$^{*}$Corresponding author.}

\author[Shuguo Zhang]{Shuguo Zhang}
\address{School of Mathematics, Sichuan University\\
Chengdu, 610064, China}
\email{zhangsg@scu.edu.cn}

\date{}

\begin{abstract}
A Hausdorff space $Y$ is called \emph{rigid} if, for every positive integer $k$, every continuous open $k$-to-one map $p:X\to Y$ from a Hausdorff space $X$ admits a partition of $X$ into $k$ open sets, each mapped homeomorphically onto $Y$. We show that the Sorgenfrey plane is not rigid, giving negative answers in $\mathsf{ZFC}$ to Velichko's questions on finite products of rigid spaces and on strongly zero-dimensional spaces. For ordinal spaces, we prove in $\mathsf{ZF}$ that every regular uncountable cardinal with the tree property is rigid. In $\mathsf{ZFC}$, we obtain a characterization for arbitrary limit ordinals: a limit ordinal $\alpha$ is rigid if and only if every coherent family $\langle f_\xi:\xi<\alpha\rangle$, where $f_\xi:\xi\to\{0,1\}$, admits a uniformization. As a consequence, a limit ordinal is rigid if and only if its cofinality is rigid. The characterization also implies that $\omega_1$ is not rigid in $\mathsf{ZFC}$, answering a question of Velichko.
\end{abstract}




\maketitle

\section{Introduction}
\label{sec:introduction}

Throughout this paper, all spaces are assumed to be Hausdorff. Unless otherwise stated, we work in $\mathsf{ZFC}$. Results proved without the Axiom of Choice are marked $\mathsf{ZF}$.

For every positive integer $k$, every continuous open $k$-to-one map is locally trivial: each point of the range has an open neighborhood whose inverse image is the disjoint union of $k$ open sets, each mapped homeomorphically onto that neighborhood~\cite[Proposition~1]{velichko2008finite}. Such a decomposition need not exist globally. Velichko~\cite{velichko2008finite} introduced the notion of a \emph{rigid space}: a space $Y$ is rigid if, for every positive integer $k$, every continuous open $k$-to-one map onto $Y$ admits a decomposition of its domain into $k$ disjoint open sets, each mapped homeomorphically onto $Y$.

Every Euclidean space and every zero-dimensional Lindel\"of space is rigid \cite[Theorems~2 and~4, Corollary~7]{velichko2008finite}. Velichko also observed that the circle is not rigid \cite[p.~S185]{velichko2008finite}. In the same paper, Velichko posed the following questions.

\begin{question}[{\cite[Question~6]{velichko2008finite}}]
\label{ques:subspace-rigid}
Is every subspace of a rigid space rigid?
\end{question}

\begin{question}[{\cite[Question~4]{velichko2008finite}}]
\label{ques:product-rigid}
Is a (finite) product of rigid spaces rigid?
\end{question}

\begin{question}[{\cite[Question~2]{velichko2008finite}}]
\label{ques:strongly-zero-dimensional}
Is every (strongly) zero-dimensional space rigid?
\end{question}

The examples mentioned above already give a negative answer to Question~\ref{ques:subspace-rigid}: the nonrigid circle is a subspace of the rigid plane $\mathbb R^2$. Our first main result (Theorem~\ref{thm:product}) shows that the Sorgenfrey plane $S\times S$ is not rigid. The Sorgenfrey line $S$ is rigid \cite[Corollary~3]{velichko2008finite}, and its square is strongly zero-dimensional \cite[Theorem~1.5]{bade1974two}. Thus this construction gives negative answers to Questions~\ref{ques:product-rigid} and~\ref{ques:strongly-zero-dimensional}.

The main part of the paper concerns the rigidity of ordinal spaces. Every successor ordinal is compact and zero-dimensional, and hence rigid in $\mathsf{ZF}$ (Proposition~\ref{prop:compact-ordinal-rigid}). Our study of limit ordinals is motivated by the following question of Velichko.

\begin{question}[{\cite[Question~1]{velichko2008finite}}]
\label{ques:omega-one-rigid}
Is the space $\omega_1$ of all countable ordinals rigid in the natural topology?
\end{question}

We give a negative answer in $\mathsf{ZFC}$ (Corollary~\ref{cor:omega-one-not-rigid}). In $\mathsf{ZF}$, however, $\omega_1$ is rigid whenever it is regular and has the tree property (Theorem~\ref{thm:tree-property-rigidity}).

For an arbitrary limit ordinal $\alpha$, we prove in $\mathsf{ZFC}$ that $\alpha$ is rigid if and only if every coherent family $\langle f_\xi:\xi<\alpha\rangle$, where $f_\xi:\xi\to\{0,1\}$, admits a uniformization (Theorem~\ref{thm:ordinal-rigidity}). This characterization clarifies several aspects of rigidity for limit ordinals:

\begin{itemize}
\item For limit ordinals, it suffices to test rigidity using continuous open $2$-to-one maps onto the ordinal.

\item Rigidity depends only on cofinality: a limit ordinal $\alpha$ is rigid if and only if $\operatorname{cf}(\alpha)$ is rigid (Corollary~\ref{cor:rigidity-cofinality}). In particular, every ordinal of cofinality $\omega$ is rigid, and every ordinal of cofinality $\omega_1$ is not rigid.

\item Known results on coherent families give further conclusions under additional set-theoretic assumptions. The $P$-ideal dichotomy implies that every ordinal of cofinality greater than $\omega_1$ is rigid, whereas under $V=L$, every ordinal whose cofinality is an uncountable successor cardinal is not rigid.
\end{itemize}

The paper is organized as follows. Section~\ref{sec:preliminaries} introduces the basic definitions and collects the preliminary results. Section~\ref{sec:preservation} studies subspaces, quotients, and proves that the Sorgenfrey plane is not rigid. Section~\ref{sec:tree-property} proves, in $\mathsf{ZF}$, that every regular uncountable cardinal with the tree property is rigid. Finally, Section~\ref{sec:combinatorial-characterization} characterizes the rigidity of limit ordinals in terms of uniformizations of coherent families and derives consequences of this characterization.

\section{Preliminaries}\label{sec:preliminaries}

Every ordinal is equipped with its order topology, and every subset of an ordinal with the subspace topology inherited from it. Our terminology and notation in general topology and set theory follow Engelking~\cite{Engelking1989} and Jech~\cite{jech2003set}, respectively.

For each positive integer $k$, we identify $k$ with the discrete space $\{0,\ldots,k-1\}$ and let $S_k$ denote the group of all permutations of $k$, also equipped with the discrete topology. We also regard $2$ as the additive group $\mathbb Z/2\mathbb Z$ and denote its addition by $+_2$ to distinguish it from ordinal addition. 

A map $p:X\to Y$ is called \emph{$k$-to-one} if $\abs{p^{-1}(\{y\})}=k$ for every $y\in Y$.

For spaces $X$ and $Y$, let $C(X,Y)$ denote the set of all continuous maps from $X$ to $Y$ and let $\operatorname{Clop}(X)$ denote the set of all clopen subsets of $X$.

\subsection{Rigid spaces}

Let $k$ be a positive integer and let $p:X\to Y$ be a continuous open $k$-to-one map.
\begin{itemize}
    \item We call $p$ \emph{elementary} if there is a tuple $(U_i)_{i<k}$ of pairwise disjoint open subsets of $X$ such that $X=\bigcup_{i<k}U_i$ and each restriction $p\!\upharpoonright U_i:U_i\to Y$ is a homeomorphism. Such a tuple is called a \emph{decomposition of $p$}. We denote the set of all decompositions of $p$ by $\mathcal D(p)$.
    \item For a nonempty subspace $A$ of $Y$, put
    \[
    X_A=p^{-1}(A),\qquad p_A=p\!\upharpoonright X_A:X_A\to A,
    \]
    and write $\mathcal D(A,p)=\mathcal D(p_A)$.
    \item For a nonempty subspace $B\subseteq A$ of $Y$ and $\mathcal U=(U_i)_{i<k}\in\mathcal D(A,p)$, put 
    \[
    \mathcal U\!\upharpoonright B =(U_i\cap X_B)_{i<k}.
    \]
    Then $\mathcal U\!\upharpoonright B\in\mathcal D(B,p)$. We say that $\mathcal U$ \emph{extends} $\mathcal V\in\mathcal D(B,p)$ if $\mathcal U\!\upharpoonright B=\mathcal V$.
\end{itemize}

For a space $Y$ and a positive integer $k$, consider the property
\begin{quote}
$R_k(Y)$:\quad For every space $X$, every continuous open $k$-to-one map $p:X\to Y$ is elementary.
\end{quote}
The space $Y$ is called \emph{rigid} if $R_k(Y)$ holds for every positive integer $k$.

The following proposition is a reformulation of \cite[Proposition~1]{velichko2008finite}, whose proof does not use the Axiom of Choice.

\begin{proposition}[$\mathsf{ZF}$]
\label{prop:local-elementary}
Let $k$ be a positive integer and let $p:X\to Y$ be a continuous open $k$-to-one map. Every $y\in Y$ has an open neighborhood $V$ such that $\mathcal D(V,p)\ne\varnothing$.
\end{proposition}
\begin{proof}
Fix $y\in Y$ and write $p^{-1}(\{y\})=\{x_i:i<k\}$. Since $X$ is Hausdorff, choose pairwise disjoint open neighborhoods $O_i$ of $x_i$ for $i<k$. Put $V=\bigcap_{i<k}p[O_i]$. Since $p$ is open, $V$ is an open neighborhood of $y$.

For each $i<k$, put $U_i=O_i\cap p^{-1}(V)$, which is open in $X$. For each $z\in V$, the fiber $p^{-1}(\{z\})$ has $k$ points and meets the $k$ pairwise disjoint sets $U_i$. Hence $p^{-1}(\{z\})\cap U_i$ is a singleton for every $i<k$. Thus $(U_i)_{i<k}$ partitions $p^{-1}(V)$, and each restriction $p\res U_i:U_i\to V$ is a continuous open bijection, hence a homeomorphism. Therefore $(U_i)_{i<k}\in\mathcal D(V,p)$.
\end{proof}

The following proposition is the compact case of \cite[Theorem~2]{velichko2008finite}. We include the proof to make explicit that the result holds in $\mathsf{ZF}$.

\begin{proposition}[$\mathsf{ZF}$]\label{prop:compact-zero-dimensional-rigid}
Every compact zero-dimensional space is rigid.
\end{proposition}

\begin{proof}
Let $Y$ be a compact zero-dimensional space, $k$ a positive integer, and $p:X\to Y$ a continuous open $k$-to-one map, where $X$ is a space. By Proposition~\ref{prop:local-elementary} and the zero-dimensionality of $Y$, every $y\in Y$ has a clopen neighborhood $V$ such that $\mathcal D(V,p)\ne\varnothing$. By compactness, finitely many such sets $V_0,\ldots,V_{m-1}$ cover $Y$.

For each $j<m$, put
\[
    W_j=V_j\setminus\bigcup_{r<j}V_r,
\]
and let $J=\{j<m:W_j\ne\varnothing\}$. Then $(W_j)_{j\in J}$ is a clopen partition of $Y$. Since $W_j\subseteq V_j$, we have $\mathcal D(W_j,p)\ne\varnothing$ for every $j\in J$. For each $j\in J$, choose $(U_{j,i})_{i<k}\in\mathcal D(W_j,p)$.

For each $i<k$, put
\[
    U_i=\bigcup_{j\in J}U_{j,i}.
\]
Then $X=\bigsqcup_{i<k}U_i$, and each $U_i$ is open in $X$. Moreover, each restriction $p\!\upharpoonright U_i:U_i\to Y$ is a continuous open bijection, hence a homeomorphism. Therefore $p$ is elementary.
\end{proof}

\subsection{Trees and the tree property}

A \emph{tree} is a set $T$ with a strict partial order $<_T$ such that, for each $x\in T$, the set $\{y\in T:y<_T x\}$ is well-ordered by $<_T$. Let $o_T(x)$ denote the order type of this set. For each ordinal $\delta$, the $\delta$th \emph{level} of $T$ is $T_\delta=\{x\in T:o_T(x)=\delta\}$. The \emph{height} of $T$ is the least ordinal $\delta$ such that $T_\delta=\varnothing$. A \emph{branch} in $T$ is a maximal linearly ordered subset
of $T$.

A regular uncountable cardinal $\kappa$ has the \emph{tree property} if every tree of height $\kappa$ whose levels have cardinality less than $\kappa$ has a branch of cardinality $\kappa$.

\subsection{Coherent families and uniformization}

For maps $f$ and $g$ with a common domain $D$, write $f=^*g$ if $\{\eta\in D:f(\eta)\ne g(\eta)\}$ is finite.

Let $\delta$ be an ordinal and $A$ a nonempty set. A family $\Phi=\langle f_\xi:\xi<\delta\rangle$, where $f_\xi:\xi\to A$ for each $\xi<\delta$, is called \emph{coherent} (see \cite[Section~3.1]{koenig2003local}) if
\[
    f_\eta\restriction\xi=^*f_\xi
\]
for every $\xi<\eta<\delta$. The family $\Phi$ is called \emph{trivial} if there exists a map $f:\delta\to A$ such that $f\restriction\xi=^*f_\xi$ for every $\xi<\delta$. Such a map is called a \emph{uniformization} of $\Phi$ and is said to \emph{uniformize} $\Phi$. We use the following notation:
\begin{quote}
$\mathrm{CU}_A(\delta)$:\quad
Every coherent family $\langle f_\xi:\xi<\delta\rangle$, where $f_\xi:\xi\to A$ for each $\xi<\delta$, is trivial.
\end{quote}

\begin{proposition}[$\mathsf{ZF}$]\label{prop:finite-valued-uniformization}
Let $\delta$ be an ordinal. If $\mathrm{CU}_2(\delta)$ holds, then $\mathrm{CU}_A(\delta)$ holds for every nonempty finite set $A$.
\end{proposition}
\begin{proof}
Assume $\mathrm{CU}_2(\delta)$. Let $A$ be a nonempty finite set and let $\Phi=\langle f_\xi:\xi<\delta\rangle$ be a coherent family, where $f_\xi:\xi\to A$ for each $\xi<\delta$. Choose a positive integer $m$ such that $\abs{A}\le 2^m$. Fix an injection $e:A\to{}^{m}2$.

For each $i<m$ and $\xi<\delta$, define $f_\xi^i:\xi\to2$ by
\[
    f_\xi^i(\eta)=e(f_\xi(\eta))(i)
\]
for every $\eta<\xi$. The coherence of $\Phi$ implies that $\Phi^i=\langle f_\xi^i:\xi<\delta\rangle$ is coherent for every $i<m$. By $\mathrm{CU}_2(\delta)$, choose a map $g_i:\delta\to2$ that uniformizes $\Phi^i$ for each $i<m$.

Define $h:\delta\to{}^{m}2$ by
\[
    h(\eta)=(g_i(\eta))_{i<m}.
\]
for every $\eta<\delta$. Fix $\xi<\delta$. For each $i<m$, we have $g_i(\eta)=f_\xi^i(\eta)=e(f_\xi(\eta))(i)$ for all but finitely many $\eta<\xi$. Since $m$ is finite, these equalities hold simultaneously for all $i<m$ outside a finite subset of $\xi$. Hence $h\restriction\xi=^*e\circ f_\xi$.

Fix $a_0\in A$ and define $g:\delta\to A$ by
\[
    g(\eta)=
    \begin{cases}
        e^{-1}(h(\eta)), & h(\eta)\in e[A],\\
        a_0, & h(\eta)\notin e[A].
    \end{cases}
\]
For each $\xi<\delta$, since $h\restriction\xi=^*e\circ f_\xi$, the definition of $g$ gives $g\restriction\xi=^*f_\xi$.
\end{proof}

\section{Preservation properties of rigidity}
\label{sec:preservation}

We examine the behavior of rigidity under passage to subspaces, quotient spaces, and finite products. 

\begin{example}\label{ex:subspaces}
Rigidity need not be inherited by open dense or closed subspaces. Recall that both $\mathbb R$ and $\mathbb R^2$ are rigid by \cite[Theorem~4 and Corollary~7]{velichko2008finite}.
\begin{itemize}
    \item The punctured plane $\mathbb R^2\setminus\{(0,0)\}$ is an open dense subspace of $\mathbb R^2$.  Consider the map
    \[
        p:\mathbb R^2\setminus\{(0,0)\}
        \longrightarrow \mathbb R^2\setminus\{(0,0)\},
        \qquad p(x,y)=(x^2-y^2,2xy).
    \]
    In polar coordinates, $p$ sends $(r,\theta)$ to $(r^2,2\theta)$. Thus $p$ is a continuous open $2$-to-one map. Since its domain is connected, $p$ is not elementary.

    \item The unit circle $\mathbb S^1=\{(x,y)\in\mathbb R^2:x^2+y^2=1\}$ is a closed subspace of $\mathbb R^2$. Velichko~\cite[p.~S185]{velichko2008finite} observed that $\mathbb S^1$ is not rigid. Indeed, the restriction
    \[
    q=p\!\upharpoonright\mathbb S^1:
    \mathbb S^1\longrightarrow\mathbb S^1
    \]
    is also a continuous open $2$-to-one map. Since $\mathbb S^1$ is connected, $q$ is not elementary.
\end{itemize}

Rigidity need not be preserved even by open quotient maps. Indeed, the map
\[
    h:\mathbb R\to\mathbb S^1,
    \qquad h(t)=(\cos 2\pi t,\sin 2\pi t),
\]
is a continuous open surjection and hence a quotient map from the rigid space $\mathbb R$ onto the nonrigid space $\mathbb S^1$.
\end{example}

The \emph{Sorgenfrey line} $S$ is the real line equipped with the topology generated by the intervals $[a,b)$ with $a<b$. The product $S\times S$ is called the \emph{Sorgenfrey plane}. Recall that the Sorgenfrey line $S$ is rigid by \cite[Corollary~3]{velichko2008finite}. The following theorem gives a negative answer to Question~\ref{ques:product-rigid}.

\begin{theorem}[$\mathsf{ZF}$]\label{thm:product}
The Sorgenfrey plane $S\times S$ is not rigid. Moreover, $R_2(S\times S)$ fails.
\end{theorem}
\begin{proof}
Put $B=S\times S$ and consider the set $X=B\times2$. Let $p:X\to B$ be the projection $p(z,i)=z$. We shall equip $X$ with a Hausdorff topology making $p$ continuous and open but not elementary.

For each $t\in\mathbb R$, put $d_t=(t,-t)\in B$. Let $D=\{d_t:t\in\mathbb R\}$ and $W=B\setminus D$. The set $D$ is closed and discrete in $B$. Let $\tau_W$ denote the subspace topology on $W$ inherited from $B$.

For each $t\in\mathbb R$ and $n<\omega$, put
\[
    U_{t,n}=[t,t+2^{-n})\times[-t,-t+2^{-n}).
\]
For each $t\in\mathbb R$, the family $\{U_{t,n}:n<\omega\}$ is a decreasing clopen neighborhood base at $d_t$, with $U_{t,0}\cap D=\{d_t\}$ and $\bigcap_{n<\omega}U_{t,n}=\{d_t\}$.

We claim that there are continuous maps $g_t:U_{t,0}\setminus\{d_t\}\to2$ for $t\in\mathbb R$ such that, for every $h\in C(W,2)$, there is an $s\in\mathbb R$ satisfying
\begin{equation}\label{eq:oscillation}
    (h+_2g_s)[U_{s,n}\setminus\{d_s\}]=2
    \qquad\text{for all }n<\omega.
\end{equation}
Indeed, the set $Q=(\mathbb Q\times\mathbb Q)\cap W$ is countable and dense in $W$, so each member of $C(W,2)$ is uniquely determined by its restriction to $Q$. Hence $\abs{C(W,2)}\le\abs{{}^{Q}2}=\abs{\mathbb R}$. Since $C(W,2)$ is nonempty, write $C(W,2)=\{h_t:t\in\mathbb R\}$. For each $t\in\mathbb R$, define $a_t:U_{t,0}\setminus\{d_t\}\to2$ by setting $a_t(z)=n\bmod2$ whenever $z\in U_{t,n}\setminus U_{t,n+1}$ for some $n<\omega$. The sets $U_{t,n}\setminus U_{t,n+1}$ form a partition of $U_{t,0}\setminus\{d_t\}$ into nonempty clopen sets, so $a_t$ is continuous. For every $n<\omega$, the map $a_t$ takes different values on $U_{t,n}\setminus U_{t,n+1}$ and $U_{t,n+1}\setminus U_{t,n+2}$. Hence $a_t[U_{t,n}\setminus\{d_t\}]=2$. For each $t\in\mathbb R$, define $g_t:U_{t,0}\setminus\{d_t\}\to2$ by $g_t(z)=h_t(z)+_2a_t(z)$. Then $g_t$ is continuous. Given $h\in C(W,2)$, choose $s\in\mathbb R$ such that $h_s=h$. For every $z\in U_{s,0}\setminus\{d_s\}$, we have $h(z)+_2g_s(z)=a_s(z)$. Thus \eqref{eq:oscillation} holds, proving the claim.

For $t\in\mathbb R$, $i<2$, and $n<\omega$, put
\[
    N(t,i,n)=\{(d_t,i)\}\cup
    \{(z,i+_2g_t(z)):z\in U_{t,n}\setminus\{d_t\}\}.
\]
Let $\mathcal B$ consist of all $V\times\{i\}$, where $V\in\tau_W$ and $i<2$, together with all $N(t,i,n)$. The family $\mathcal B$ covers $X$. By continuity of $g_t$, whenever $(z,j)\in N(t,i,n)$ with $z\in W$, there is a $V\in\tau_W$ containing $z$ such that $V\times\{j\}\subseteq N(t,i,n)$. Moreover, $N(t,i,n)\cap(D\times2)=\{(d_t,i)\}$, and the sets $N(t,i,n)$ decrease with $n$. Hence $\mathcal B$ is a base; equip $X$ with the topology it generates.

For any $V\in\tau_W$ and $i<2$,
\[
    p[V\times\{i\}]=V,\qquad p^{-1}(V)=V\times2.
\]
For any $t\in\mathbb R$, $i<2$, and $n<\omega$, we also have
\[
    p[N(t,i,n)]=U_{t,n},\qquad
    p^{-1}(U_{t,n})=N(t,0,n)\sqcup N(t,1,n).
\]
As $W$ is open in $B$ and, for each $t\in\mathbb R$, $\{U_{t,n}:n<\omega\}$ is a neighborhood base at $d_t$, these identities show that $p$ is open and continuous. To show that $X$ is Hausdorff, it suffices to separate  $(z,0)$ and $(z,1)$ for each $z\in B$ since $B$ is Hausdorff and $p$ is continuous. If $z\in W$, they are separated by $W\times\{0\}$ and $W\times\{1\}$; if $z=d_t\in D$, by $N(t,0,0)$ and $N(t,1,0)$. Thus $X$ is Hausdorff.

Suppose that $p$ is elementary, and let $\mathcal U=(U_0,U_1)\in\mathcal D(p)$. For each $z\in W$, let $h(z)$ be the unique $i<2$ such that $(z,i)\in U_0$. Then $h\in C(W,2)$, since
\[
    h^{-1}(\{i\})=p[U_0\cap(W\times\{i\})]
\]
is open in $W$ for each $i<2$. Choose $s\in\mathbb R$ satisfying \eqref{eq:oscillation} for this $h$, and let $i_0<2$ be the unique index such that $(d_s,i_0)\in U_0$. Since $U_0$ is open, $N(s,i_0,n_0)\subseteq U_0$ for some $n_0<\omega$. It follows that
\[
    h(z)=i_0+_2g_s(z)
    \qquad\text{for every }z\in U_{s,n_0}\setminus\{d_s\},
\]
and hence $(h+_2g_s)[U_{s,n_0}\setminus\{d_s\}]=\{i_0\}$, contrary to \eqref{eq:oscillation}. Therefore $p$ is not elementary, and $R_2(B)$ fails.
\end{proof}

Since the Sorgenfrey plane $S\times S$ is strongly zero-dimensional \cite[Theorem~1.5]{bade1974two}, Theorem~\ref{thm:product} also gives the following negative answer to Question~\ref{ques:strongly-zero-dimensional}.

\begin{corollary}\label{cor:strongly-zero-dimensional}
Not every strongly zero-dimensional space is rigid.
\end{corollary}

\section{The tree property and rigidity of ordinals}
\label{sec:tree-property}

In this section, we prove that every regular uncountable cardinal with the tree property is rigid.

Every successor ordinal is compact and zero-dimensional. Thus Proposition~\ref{prop:compact-zero-dimensional-rigid} yields the following.

\begin{proposition}[$\mathsf{ZF}$]\label{prop:compact-ordinal-rigid}
Every successor ordinal is rigid.
\end{proposition}

We first bound the number of decompositions over a successor ordinal.

\begin{lemma}[$\mathsf{ZF}$]\label{lem:ordinal-decompositions}
Let $\delta$ be an ordinal, $k$ a positive integer, and $p:X\to\delta+1$ a continuous open $k$-to-one map, where $X$ is a space. Then $\mathcal D(p)$ is nonempty and
\[
    \abs{\mathcal D(p)}\le\max\{\aleph_0,\abs{\delta}\}.
\]
\end{lemma}
\begin{proof}
By Proposition~\ref{prop:compact-ordinal-rigid}, $\mathcal D(p)$ is nonempty. Fix a decomposition $\mathcal U=(U_i)_{i<k}\in\mathcal D(p)$ and put $\mu=\max\{\aleph_0,\abs{\delta}\}$. For each $\xi\le\delta$ and $i<k$, let $x_i^\xi$ be the unique point of $U_i\cap p^{-1}(\{\xi\})$. Then $p^{-1}(\{\xi\})=\{x_i^\xi:i<k\}$ for every $\xi\le\delta$. Thus $\mathcal U$ gives a fixed labeling of the points in each fiber of $p$.

Now let $\mathcal V=(V_j)_{j<k}\in\mathcal D(p)$. Define $\sigma_{\mathcal V}:\delta+1\to S_k$ by
\[
    \sigma_{\mathcal V}(\xi)(i)=j
    \quad\Longleftrightarrow\quad
    x_i^\xi\in V_j
\]
for every $\xi\le\delta$ and $i,j<k$. Since each $V_j$ meets every fiber of $p$ in exactly one point, $\sigma_{\mathcal V}(\xi)\in S_k$ for every $\xi\le\delta$.

The assignment $\mathcal V\mapsto\sigma_{\mathcal V}$
is injective. Indeed, for every $j<k$,
\[
    V_j=\bigcup_{i<k}
    \{x\in U_i:\sigma_{\mathcal V}(p(x))(i)=j\}.
\]
We next show that $\sigma_{\mathcal V}$ is continuous. For every $\tau\in S_k$,
\[
    \sigma_{\mathcal V}^{-1}(\{\tau\}) =\bigcap_{i<k}p[U_i\cap V_{\tau(i)}],
\]
which is open because $p$ is open. Since $S_k$ is discrete, $\sigma_{\mathcal V}$ is continuous. Thus $\abs{\mathcal D(p)}\le\abs{C(\delta+1,S_k)}$.

It remains to bound $\abs{C(\delta+1,S_k)}$. Put $\mathcal B =\{[0,\eta]:\eta\le\delta\} \cup \{(\xi,\eta]:\xi<\eta\le\delta\}$. Then $\mathcal B$ is a clopen base for $\delta+1$. Since its members can be encoded by elements of $2\times(\delta+1)^2$, the set $[\mathcal B]^{<\omega}$ is well-orderable and has cardinality at most $\mu$. 

Every member of $\Clop(\delta+1)$ is compact and hence a finite union of members of $\mathcal B$. Fix a well-ordering of $[\mathcal B]^{<\omega}$ and assign to each $C\in\Clop(\delta+1)$ the least finite subfamily whose union is $C$. This gives an injection from $\Clop(\delta+1)$ into $[\mathcal B]^{<\omega}$, so $\abs{\Clop(\delta+1)}\le\mu$.

Since $S_k$ is finite and discrete, each continuous map from $\delta+1$ to $S_k$ is determined by its clopen fibers. Consequently,
\[
    \abs{C(\delta+1,S_k)}
    \le\abs{\Clop(\delta+1)}^{k!}
    \le\mu.
\]
Therefore $\abs{\mathcal D(p)}\le\mu$.
\end{proof}

\begin{theorem}[$\mathsf{ZF}$]\label{thm:tree-property-rigidity}
Let $\kappa$ be a regular uncountable cardinal. If $\kappa$ has the tree property, then $\kappa$ is rigid.
\end{theorem}
\begin{proof}
Let $k$ be a positive integer and $p:X\to\kappa$ a continuous open $k$-to-one map, where $X$ is a space. For each $\delta<\kappa$, the initial segment $\delta+1$ is open in $\kappa$, so $X_{\delta+1}$ is open in $X$ and $p_{\delta+1}:X_{\delta+1}\to\delta+1$ is continuous, open and $k$-to-one.

We organize the decompositions of the $p_{\delta+1}$ into a tree. Put
\[
    T=\bigcup_{\delta<\kappa}\mathcal D(p_{\delta+1}).
\]
Since $p$ is surjective, the sets $X_{\delta+1}$ are pairwise distinct. Thus each member of $T$ belongs to a unique $\mathcal D(p_{\delta+1})$. For $\xi,\eta<\kappa$, $\mathcal U\in\mathcal D(p_{\xi+1})$, and $\mathcal V\in\mathcal D(p_{\eta+1})$, define
\[
    \mathcal U<_T\mathcal V \quad\Longleftrightarrow\quad \xi<\eta  \text{ and } \mathcal U=\mathcal V\res(\xi+1).
\]
For each $\eta<\kappa$ and $\mathcal V\in\mathcal D(p_{\eta+1})$, its predecessors are precisely $\mathcal V\res(\xi+1)$ for $\xi<\eta$, and their order type under $<_T$ is $\eta$. Thus $(T,<_T)$ is a tree with $T_\delta=\mathcal D(p_{\delta+1})$ for every $\delta<\kappa$. By Lemma~\ref{lem:ordinal-decompositions}, for every $\delta<\kappa$, the level $T_\delta$ is nonempty and
\[
    \abs{T_\delta} \le\max\{\aleph_0,\abs{\delta}\}  <\kappa.
\]
Hence $T$ has height $\kappa$. 

By the tree property, $T$ has a branch $b$ of cardinality $\kappa$. The branch $b$ meets each level in exactly one node. For each $\delta<\kappa$, let $\mathcal U^\delta=(U_i^\delta)_{i<k}$ be the unique member of $b\cap T_\delta$. Since $b$ is a chain, the definition of $<_T$ gives
\[
    U_i^\xi=U_i^\eta\cap X_{\xi+1} 
\]
for every $i<k$ and $\xi<\eta<\kappa$.

For each $i<k$, put
\[
    U_i=\bigcup_{\delta<\kappa}U_i^\delta.
\]
Each $U_i$ is open in $X$, and the compatibility of the decompositions gives $U_i\cap X_{\delta+1}=U_i^\delta$ for every $i<k$ and $\delta<\kappa$. Since $X=\bigcup_{\delta<\kappa}X_{\delta+1}$ and each $\mathcal U^\delta$ is a decomposition of $p_{\delta+1}$, the family $\{U_i:i<k\}$ partitions $X$, and each fiber of $p$ meets each $U_i$ in exactly one point. Thus each restriction $p\res U_i:U_i\to\kappa$ is a continuous open bijection, hence a homeomorphism. Therefore $(U_i)_{i<k}\in\mathcal D(p)$. Since $k$ and $p$ are arbitrary, $\kappa$ is rigid.
\end{proof}

In $\mathsf{ZF}$, the assumptions $\mathsf{AD}$ and
$V=L(\mathbb R)$ imply that $\omega_1$ is regular
and has the tree property; see \cite[Theorem~3.4]{ChanJackson2019}. Thus Theorem~\ref{thm:tree-property-rigidity} yields the following.

\begin{corollary}[$\mathsf{ZF}$]\label{cor:omega-one-rigid-determinacy}
Assume $\mathsf{AD}$ and $V=L(\mathbb R)$. Then $\omega_1$ is rigid.
\end{corollary}

\section{A combinatorial characterization of rigidity}
\label{sec:combinatorial-characterization}

We characterize the rigidity of limit ordinals in terms of coherent families of maps.

\begin{theorem}\label{thm:ordinal-rigidity}
Let $\alpha$ be a limit ordinal. Then the following statements are equivalent:
\begin{enumerate}[label=\upshape(\roman*), leftmargin=*, widest=iii]
    \item\label{it:ordinal-rigid} $\alpha$ is rigid;
    \item\label{it:ordinal-r-two} $R_2(\alpha)$ holds;
    \item\label{it:ordinal-cu-two} $\mathrm{CU}_2(\alpha)$ holds;
    \item\label{it:ordinal-cu-finite} $\mathrm{CU}_A(\alpha)$ holds for every nonempty finite set $A$.
\end{enumerate}
\end{theorem}

The proof uses two constructions. Local decompositions of a continuous open $k$-to-one map onto $\alpha$ yield a coherent family that is trivial if and only if the map is elementary. Conversely, given a coherent family of maps into $2$ indexed by $\alpha$, we construct a continuous open $2$-to-one map onto $\alpha$ that is elementary if and only if the family is trivial.

\begin{lemma}[$\mathsf{ZF}$]\label{lem:section5-coding}
Let $\alpha$ be a limit ordinal, $k$ a positive integer, and $p:X\to\alpha$ a continuous open $k$-to-one map, where $X$ is a space. Suppose that a family $\langle\mathcal U^\xi:\xi<\alpha\rangle$ is given, where $\mathcal U^\xi=(U_i^\xi)_{i<k}\in\mathcal D(\xi+1,p)$ for every $\xi<\alpha$. For each $\xi<\alpha$, define $b_\xi:k\to p^{-1}(\{\xi\})$ by letting $b_\xi(i)$ be the unique point of $U_i^\xi\cap p^{-1}(\{\xi\})$ for every $i<k$. 

For $\xi<\alpha$ and $\mathcal U=(U_i)_{i<k}\in\mathcal D(\xi+1,p)$, define the associated map $f_{\mathcal U}:\xi\to S_k$ by
\begin{equation}\label{eq:local-code}
    f_{\mathcal U}(\eta)(r)=s \quad\Longleftrightarrow\quad b_\eta(r),b_{\eta+1}(s)\in U_i \text{ for some }i<k,
\end{equation}
where $\eta<\xi$ and $r,s<k$.

Put $\Phi=\langle f_{\mathcal U^\xi}:\xi<\alpha\rangle$.
Then the following statements hold:
\begin{enumerate}[label=\upshape(\roman*), leftmargin=*, widest=ii]
    \item\label{it:coding-coherent} $\Phi$ is coherent.
    \item\label{it:coding-elementary} $p$ is elementary if and only if $\Phi$ is trivial.
\end{enumerate}
\end{lemma}

\begin{proof}
Each $b_\xi$ is a bijection determined by $\mathcal U^\xi$. Since each member of a decomposition meets every fiber in exactly one point, \eqref{eq:local-code} defines a permutation of $k$ at each $\eta<\xi$.
\begin{claim}\label{claim:local-code-coherence}
For any $\xi\le\eta<\alpha$ and any decompositions $\mathcal U=(U_i)_{i<k}\in\mathcal D(\xi+1,p)$ and $\mathcal V=(V_i)_{i<k}\in\mathcal D(\eta+1,p)$, we have $f_{\mathcal V}\restriction\xi=^*f_{\mathcal U}$.
\end{claim}

\begin{proof}[Proof of Claim~\ref{claim:local-code-coherence}]
Suppose that $F=\{\zeta<\xi:f_{\mathcal U}(\zeta)\ne f_{\mathcal V}(\zeta)\}$ is infinite. Choose a limit ordinal $\lambda\le\xi$ such that $F\cap\lambda$ is cofinal in $\lambda$. There is a unique $\pi\in S_k$ such that
\[
    \lambda\in p[U_i\cap V_{\pi(i)}]
\]
for every $i<k$. Since $X_{\xi+1}$ and $X_{\eta+1}$ are open in $X$, the sets $U_i$ and $V_{\pi(i)}$ are open in $X$ for every $i<k$. As $p$ is open, the set $O=\bigcap_{i<k}p[U_i\cap V_{\pi(i)}]$ is an open neighborhood of $\lambda$. Choose $\rho<\lambda$ such that $(\rho,\lambda]\subseteq O$. 

For each $i<k$, the sets $U_i$ and $V_{\pi(i)}$ meet each fiber over $O$ in the same point. If $\rho<\zeta<\lambda$, then $\zeta,\zeta+1\in(\rho,\lambda]\subseteq O$, so \eqref{eq:local-code} gives $f_{\mathcal U}(\zeta)=f_{\mathcal V}(\zeta)$. This contradicts the cofinality of $F\cap\lambda$ in $\lambda$.
\end{proof}

Applying Claim~\ref{claim:local-code-coherence} to $\mathcal U^\xi$ and $\mathcal U^\eta$ for all $\xi<\eta<\alpha$ shows that $\Phi$ is coherent, proving \ref{it:coding-coherent}.

For \ref{it:coding-elementary}, suppose first that $p$ is elementary, and fix $\mathcal U=(U_i)_{i<k}\in\mathcal D(p)$. Define $f_{\mathcal U}:\alpha\to S_k$ by the same rule as in~\eqref{eq:local-code}. For each $\xi<\alpha$, the restriction $\mathcal U\restriction(\xi+1)=(U_i\cap X_{\xi+1})_{i<k}$ belongs to $\mathcal D(\xi+1,p)$, and $f_{\mathcal U\restriction(\xi+1)} =f_{\mathcal U}\restriction\xi$. Applying Claim~\ref{claim:local-code-coherence} to $\mathcal U\restriction(\xi+1)$ and $\mathcal U^\xi$ gives $f_{\mathcal U\restriction(\xi+1)}=^*f_{\mathcal U^\xi}$, which implies $f_{\mathcal U}\restriction\xi=^*f_{\mathcal U^\xi}$. Thus $f_{\mathcal U}$ uniformizes $\Phi$.

Conversely, suppose that $\Phi$ is trivial, and fix a map $f:\alpha\to S_k$ which uniformizes $\Phi$.

\begin{claim}\label{claim:local-code-realization}
For each $\xi<\alpha$, there is a unique decomposition
$\mathcal V^\xi=(V_i^\xi)_{i<k}\in\mathcal D(\xi+1,p)$ such that
\begin{itemize}
    \item $f_{\mathcal V^\xi}=f\restriction\xi$; and
    \item $b_0(i)\in V_i^\xi$ for every $i<k$.
\end{itemize}
\end{claim}

Assuming Claim~\ref{claim:local-code-realization}, we complete the proof. For each $\xi<\alpha$, let $\mathcal V^\xi=(V_i^\xi)_{i<k}$ be the unique decomposition given by Claim~\ref{claim:local-code-realization}. For $\xi<\eta<\alpha$, the restriction $\mathcal V^\eta\restriction(\xi+1)$ belongs to $\mathcal D(\xi+1,p)$, its associated map is $f\restriction\xi$, and its $i$th member contains $b_0(i)$ for every $i<k$. By uniqueness,
\[
    \mathcal V^\eta\restriction(\xi+1)=\mathcal V^\xi.
\]
For each $i<k$, put $V_i=\bigcup_{\xi<\alpha}V_i^\xi$. By compatibility, the family $\{V_i:i<k\}$ is an open partition of $X$, and each meets every fiber of $p$ in exactly one point. Since $p$ is continuous and open, each restriction $p\res V_i:V_i\to\alpha$ is a continuous open bijection, hence a homeomorphism. We obtain $(V_i)_{i<k}\in\mathcal D(p)$, so $p$ is elementary.

It remains to prove Claim~\ref{claim:local-code-realization}.

\begin{proof}[Proof of Claim~\ref{claim:local-code-realization}]
Fix $\xi<\alpha$. Since $f$ uniformizes $\Phi$, the set $F=\{\eta<\xi:f_{\mathcal U^\xi}(\eta)\ne f(\eta)\}$ is finite. Write $F=\{\eta_j:j<m\}$ in increasing order, where $m<\omega$.
If $m>0$, partition $\xi+1$ into the clopen intervals
\[
    I_0=[0,\eta_0],\qquad
    I_j=(\eta_{j-1},\eta_j]\quad\text{for }0<j<m,\qquad
    I_m=(\eta_{m-1},\xi].
\]
If $m=0$, put $I_0=\xi+1$.

We reindex $\mathcal U^\xi\restriction I_j$ successively for $j\le m$ so that the associated map of the resulting decomposition agrees with $f$ on $F$ and with $f_{\mathcal U^\xi}$ on $\xi\setminus F$. For each $i<k$, put
\[
    W_{0,i}=U_s^\xi\cap X_{I_0},
\]
where $U_s^\xi$ is the unique member of $\mathcal U^\xi$ containing $b_0(i)$. Suppose that $j<m$ and $(W_{j,i})_{i<k}$ has been defined. For each $i<k$, let $r<k$ be the unique index such that $b_{\eta_j}(r)\in W_{j,i}$, and put
\[
    W_{j+1,i}=U_{s'}^\xi\cap X_{I_{j+1}},
\]
where $U_{s'}^\xi$ is the unique member of $\mathcal U^\xi$ containing $b_{\eta_j+1}(f(\eta_j)(r))$. For each $i<k$, put
\[
    V_i^\xi=\bigcup_{j\le m}W_{j,i}.
\]
Since each $(W_{j,i})_{i<k}$ is a reindexing of $\mathcal U^\xi\restriction I_j$ and the intervals $I_j$ form a clopen partition of $\xi+1$, we have $\mathcal V^\xi=(V_i^\xi)_{i<k}\in\mathcal D(\xi+1,p)$.

For $\zeta\in\xi\setminus F$, both $\zeta$ and $\zeta+1$ lie in the same interval, so $f_{\mathcal V^\xi}(\zeta)=f_{\mathcal U^\xi}(\zeta)=f(\zeta)$. For each $j<m$ and all $i,r<k$, the construction gives
\[
    b_{\eta_j}(r)\in W_{j,i} \quad\Longleftrightarrow\quad  b_{\eta_j+1}(f(\eta_j)(r))\in W_{j+1,i},
\]
and hence $f_{\mathcal V^\xi}(\eta_j)=f(\eta_j)$. Thus $f_{\mathcal V^\xi}=f\restriction\xi$. Moreover, $b_0(i)\in W_{0,i}\subseteq V_i^\xi$ for every $i<k$.

For uniqueness, let $\mathcal V'=(V_i')_{i<k}\in\mathcal D(\xi+1,p)$ satisfy the same conditions, and put $A_i=p[V_i^\xi\cap V_i']$ for every $i<k$. To prove $\mathcal V'=\mathcal V^\xi$, it suffices to show $A_i=\xi+1$ for every $i<k$. Fix $i<k$. We verify that
\begin{itemize}
    \item $A_i$ is closed in $\xi+1$;
    \item $0\in A_i$; and
    \item $\zeta\in A_i$ implies $\zeta+1\in A_i$
    for every $\zeta<\xi$.
\end{itemize}
Since $V_i'$ is closed in $X_{\xi+1}$, the intersection $V_i^\xi\cap V_i'$ is closed in $V_i^\xi$. The restriction $p\res V_i^\xi:V_i^\xi\to\xi+1$ is a homeomorphism, so $A_i=p[V_i^\xi\cap V_i']$ is closed in $\xi+1$. Also, $b_0(i)\in V_i^\xi\cap V_i'$ gives $0\in A_i$. If $\zeta\in A_i$ and $\zeta<\xi$, choose $r<k$ such that $b_\zeta(r)\in V_i^\xi\cap V_i'$. Since $f_{\mathcal V^\xi}(\zeta)=f_{\mathcal V'}(\zeta)=f(\zeta)$,
we have
\[
    b_{\zeta+1}(f(\zeta)(r))\in V_i^\xi\cap V_i',
\]
so $\zeta+1\in A_i$. These properties imply that $A_i=\xi+1$. Since $i<k$ is arbitrary, $\mathcal V'=\mathcal V^\xi$.
\end{proof}
\end{proof}

\begin{lemma}[$\mathsf{ZF}$]\label{lem:binary-realization}
Let $\alpha$ be a limit ordinal and $\Phi=\langle f_\xi:\xi<\alpha\rangle$ a coherent family, where $f_\xi:\xi\to2$ for every $\xi<\alpha$. There is a Hausdorff topology on $X=\alpha\times2$ such that the projection $p:X\to\alpha$, defined by $p(\xi,i)=\xi$, is continuous and open, and $p$ is elementary if and only if $\Phi$ is trivial.
\end{lemma}
\begin{proof}
We first use $\Phi$ to construct two copies of each successor initial segment of $\alpha$, which will determine a topology on $X$ in which each copy is open. We then compare the family associated with these local decompositions with $\Phi$.

For $\xi\le\eta<\alpha$, put
\[
    D_{\xi\eta}
    =\{\zeta<\xi:f_\xi(\zeta)\ne f_\eta(\zeta)\},
    \qquad
    \varepsilon_{\xi\eta}=\abs{D_{\xi\eta}}\bmod2.
\]
Each $D_{\xi\eta}$ is finite by coherence. For $\eta<\alpha$ and $i<2$, define
\[
    U_i^\eta
    =\{(\xi,i+_2\varepsilon_{\xi\eta}):\xi\le\eta\}.
\]
Then $U_0^\eta$ and $U_1^\eta$ partition $X_{\eta+1}$, and each restriction $p\restriction U_i^\eta:U_i^\eta\to\eta+1$ is bijective. Declare $G\subseteq X$ open precisely when, for every $\eta<\alpha$ and $i<2$,
\[
    p[G\cap U_i^\eta]\text{ is open in }\eta+1.
\]
Since each $p\restriction U_i^\eta$ is bijective, this defines a topology on $X$. For an open set $V\subseteq\alpha$, we have $p[p^{-1}(V)\cap U_i^\eta]=V\cap(\eta+1)$ for every $\eta<\alpha$ and $i<2$. Thus $p$ is continuous. For an open set $G\subseteq X$,
\[
    p[G]=\bigcup_{\eta<\alpha}\bigcup_{i<2}p[G\cap U_i^\eta].
\]
Each set in this union is open in $\eta+1$ and hence in $\alpha$. Thus $p$ is open.

Now we verify that $\mathcal U^\eta=(U_0^\eta,U_1^\eta)$ belongs to $\mathcal D(\eta+1,p)$ for every $\eta<\alpha$ so that we can apply Lemma~\ref{lem:section5-coding}. We need only show that each $U_i^\eta$ is open in $X$. By the definition of the topology, it suffices to verify that $p[U_i^\xi\cap U_j^\eta]$ is open in $\alpha$ for $\xi\le\eta<\alpha$ and $i,j<2$. Fix such $\xi,\eta,i,j$. Then
\[
\begin{aligned}
    p[U_i^\xi\cap U_j^\eta]
    &=\{\zeta\le\xi:
        i+_2\varepsilon_{\zeta\xi}
        =j+_2\varepsilon_{\zeta\eta}\}\\
    &=\{\zeta\le\xi:
        \varepsilon_{\zeta\xi}+_2\varepsilon_{\zeta\eta}
        =i+_2j\}.
\end{aligned}
\]
Since the maps in $\Phi$ take values in $2$, we have $D_{\zeta\xi}\mathbin{\triangle}D_{\zeta\eta} =D_{\xi\eta}\cap\zeta$ for every $\zeta\le\xi$. It follows that
\begin{equation}\label{eq:binary-overlap-parity}
    \varepsilon_{\zeta\xi}+_2\varepsilon_{\zeta\eta}
    =\abs{D_{\xi\eta}\cap\zeta}\bmod2
    \qquad\text{for every }\zeta\le\xi.
\end{equation}
Let $\lambda\in p[U_i^\xi\cap U_j^\eta]$ be a limit ordinal. Since $D_{\xi\eta}$ is finite, there is a $\rho<\lambda$ such that $D_{\xi\eta}\cap\zeta=D_{\xi\eta}\cap\lambda$ whenever $\rho<\zeta\le\lambda$. Combining this equality with \eqref{eq:binary-overlap-parity},
we obtain
\[
    \varepsilon_{\zeta\xi}+_2\varepsilon_{\zeta\eta}
    =\varepsilon_{\lambda\xi}+_2\varepsilon_{\lambda\eta}=i+_2j
\]
for every $\zeta\in(\rho,\lambda]$. So $(\rho,\lambda]\subseteq p[U_i^\xi\cap U_j^\eta]$. The point $0$ and all successor ordinals in $\xi+1$ are isolated. Thus $p[U_i^\xi\cap U_j^\eta]$ is open in $\xi+1$, and hence in $\alpha$. It follows that each $U_i^\eta$ is open in $X$, and $\mathcal U^ \eta\in\mathcal D(\eta+1,p)$ for every $\eta<\alpha$.

Since $\alpha$ is Hausdorff and $p$ is continuous, points in distinct fibers can be separated by inverse images of disjoint open subsets of $\alpha$. For each $\eta<\alpha$, the points $(\eta,0)$ and $(\eta,1)$ are separated by the disjoint open sets $U_0^\eta$ and $U_1^\eta$. Hence $X$ is Hausdorff.

It remains to relate the family associated with $\langle\mathcal U^\eta:\eta<\alpha\rangle$ to $\Phi$. Since $\varepsilon_{\xi\xi}=0$, we have $U_i^\xi\cap p^{-1}(\{\xi\})=\{(\xi,i)\}$. Thus the bijections in Lemma~\ref{lem:section5-coding} are given by $b_\xi(r)=(\xi,r)$ for $\xi<\alpha$ and $r<2$. Identify each $t\in2$ with the permutation $r\mapsto r+_2t$ in $S_2$. For each $\eta<\alpha$, let $h_\eta:\eta\to2$ be the map associated with $\mathcal U^\eta$ by \eqref{eq:local-code}. For every $\xi<\eta<\alpha$, we have
\[
    h_\eta(\xi)
    =\varepsilon_{\xi\eta}+_2\varepsilon_{\xi+1,\eta}.
\]
We claim that, for any $\xi\le\eta<\alpha$, the maps $h_\xi+_2f_\xi$ and $h_\eta+_2f_\eta$ agree on $\xi$. Fix $\zeta<\xi\le\eta<\alpha$. By the formula for $h_\xi$ and $h_\eta$, applying \eqref{eq:binary-overlap-parity} at $\zeta$ and $\zeta+1$ gives
\[
\begin{aligned}
    h_\xi(\zeta)+_2h_\eta(\zeta)
    &=(\varepsilon_{\zeta\xi}+_2\varepsilon_{\zeta+1,\xi})+_2(\varepsilon_{\zeta\eta}+_2\varepsilon_{\zeta+1,\eta})\\
    &=(\varepsilon_{\zeta\xi}+_2\varepsilon_{\zeta\eta})+_2(\varepsilon_{\zeta+1,\xi}+_2\varepsilon_{\zeta+1,\eta})\\
    &=\bigl(\abs{D_{\xi\eta}\cap\zeta}\bmod2\bigr)
      +_2\bigl(\abs{D_{\xi\eta}\cap(\zeta+1)}\bmod2\bigr)\\
    &=\abs{D_{\xi\eta}\cap\{\zeta\}}\bmod2\\
    &=f_\xi(\zeta)+_2f_\eta(\zeta).
\end{aligned}
\]
So $h_\xi(\zeta)+_2f_\xi(\zeta)=h_\eta(\zeta)+_2f_\eta(\zeta)$. Since $\alpha$ is a limit ordinal, the union
\[
    a=\bigcup_{\eta<\alpha}(h_\eta+_2f_\eta)
\] 
is a map from $\alpha$ to $2$ and satisfies $a\restriction\eta=h_\eta+_2f_\eta$ for every $\eta<\alpha$. A map $g:\alpha\to2$ uniformizes $\langle h_\eta:\eta<\alpha\rangle$ if and only if $g+_2a$ uniformizes $\Phi$. By Lemma~\ref{lem:section5-coding}, $p$ is elementary if and only if $\Phi$ is trivial.
\end{proof}

\begin{proof}[Proof of Theorem~\ref{thm:ordinal-rigidity}]
The implication \ref{it:ordinal-rigid}$\Rightarrow$\ref{it:ordinal-r-two} follows from the definition of rigidity.

For \ref{it:ordinal-r-two}$\Rightarrow$\ref{it:ordinal-cu-two}, let $\Phi=\langle f_\xi:\xi<\alpha\rangle$ be a coherent family, where $f_\xi:\xi\to2$ for every $\xi<\alpha$. By Lemma~\ref{lem:binary-realization}, there are a space $X$ and a continuous open $2$-to-one map $p:X\to\alpha$ such that $p$ is elementary if and only if $\Phi$ is trivial. By \ref{it:ordinal-r-two}, $p$ is elementary, so $\Phi$ is trivial.

The implication \ref{it:ordinal-cu-two}$\Rightarrow$\ref{it:ordinal-cu-finite} follows from Proposition~\ref{prop:finite-valued-uniformization}.

Finally, assume \ref{it:ordinal-cu-finite}. Let $k$ be a positive integer and let $p:X\to\alpha$ be a continuous open $k$-to-one map, where $X$ is a space. By Proposition~\ref{prop:compact-ordinal-rigid}, $\mathcal D(\xi+1,p)$ is nonempty for every $\xi<\alpha$. Using the Axiom of Choice, choose a family $\langle\mathcal U^\xi:\xi<\alpha\rangle$ with $\mathcal U^\xi\in\mathcal D(\xi+1,p)$ for every $\xi<\alpha$. By Lemma~\ref{lem:section5-coding}, the family $\Phi=\langle f_{\mathcal U^\xi}:\xi<\alpha\rangle$ is coherent, and $p$ is elementary if and only if $\Phi$ is trivial. Since $S_k$ is finite and nonempty, \ref{it:ordinal-cu-finite} implies that $\Phi$ is trivial. Hence $p$ is elementary. Since $k$ and $p$ are arbitrary, $\alpha$ is rigid.
\end{proof}

By \cite[Lemma~7.5]{bergfalk2021}, there is a nontrivial coherent family $\langle\varphi_\eta:\eta\in E\rangle$, where $E\subseteq\omega_1$ is club and $\varphi_\eta:\eta\to2$. For each $\xi<\omega_1$, put $f_\xi=\varphi_{\min(E\setminus\xi)}\restriction\xi$. Then $\langle f_\xi:\xi<\omega_1\rangle$ is coherent, and any uniformization of this family would uniformize the original family. Hence $\mathrm{CU}_2(\omega_1)$ fails. Theorem~\ref{thm:ordinal-rigidity} therefore yields the following corollary.

\begin{corollary}\label{cor:omega-one-not-rigid}
The space $\omega_1$ is not rigid.
\end{corollary}

The following proposition reduces the study of $\mathrm{CU}_2(\alpha)$ to the case of regular cardinals.

\begin{proposition}[$\mathsf{ZF}$]\label{prop:cu-cofinality}
Let $\alpha$ be a limit ordinal. Then $\mathrm{CU}_2(\alpha)$ holds if and only if $\mathrm{CU}_2(\operatorname{cf}(\alpha))$ holds.
\end{proposition}
\begin{proof}
Put $\kappa=\operatorname{cf}(\alpha)$, and fix a strictly increasing cofinal sequence $\langle c_\xi:\xi<\kappa\rangle$ in $\alpha$.

Suppose that $\mathrm{CU}_2(\alpha)$ holds, and let $\langle f_\xi:\xi<\kappa\rangle$ be a coherent family, where $f_\xi:\xi\to2$ for every $\xi<\kappa$. For each $\xi<\kappa$, define $g_{c_\xi}:c_\xi\to2$ by
\[
    g_{c_\xi}(c_\eta)=f_\xi(\eta)
\]
for every $\eta<\xi$, and set $g_{c_\xi}(\delta)=0$
at all remaining points of $c_\xi$. For each $\delta\in\alpha\setminus\{c_\eta:\eta<\kappa\}$, put $g_\delta=g_{c_\xi}\restriction\delta$, where $\xi$ is the least index such that $\delta\le c_\xi$. By the coherence of $\langle f_\xi:\xi<\kappa\rangle$, we have $g_{c_\eta}\restriction c_\xi=^*g_{c_\xi}$ whenever $\xi<\eta<\kappa$. It follows from the definition of the maps $g_\delta$ that $\langle g_\delta:\delta<\alpha\rangle$ is coherent. By $\mathrm{CU}_2(\alpha)$, choose a map $g:\alpha\to2$ that uniformizes $\langle g_\delta:\delta<\alpha\rangle$. Define $f:\kappa\to2$ by $f(\xi)=g(c_\xi)$. For each $\xi<\kappa$, the relation $g\restriction c_\xi=^*g_{c_\xi}$ and the injectivity of $\eta\mapsto c_\eta$ give $f\restriction\xi=^*f_\xi$. Thus $\mathrm{CU}_2(\kappa)$ holds.

Conversely, suppose that $\mathrm{CU}_2(\kappa)$ holds, and let $\langle g_\delta:\delta<\alpha\rangle$ be a coherent family, where $g_\delta:\delta\to2$ for every $\delta<\alpha$. By cofinality and coherence, it suffices to find a map $g:\alpha\to2$ such that
\[
    g\restriction c_\eta=^*g_{c_\eta} 
\]
for every $\eta<\kappa$.

\noindent\textbf{Case 1.} $\kappa=\omega$.
Define $g:\alpha\to2$ by
\[
\begin{aligned}
    g\restriction c_0&=g_{c_0},\\
    g\restriction[c_\xi,c_{\xi+1})
    &=g_{c_{\xi+1}}\restriction[c_\xi,c_{\xi+1})
    \qquad\text{for every }\xi<\omega.
\end{aligned}
\]
For each $\eta<\omega$, $g\restriction c_\eta$ is obtained by pasting finitely many restrictions of $g_{c_\xi}$ with $\xi\le\eta$. Thus coherence gives $g\restriction c_\eta=^*g_{c_\eta}$.

\noindent\textbf{Case 2.} $\kappa>\omega$.
Put
\[
    B=\{c_\xi:\xi<\kappa\}\cup
      \bigcup_{\xi<\eta<\kappa}
      \{\delta<c_\xi:
        g_{c_\xi}(\delta)\ne g_{c_\eta}(\delta)\}.
\]
Then $\abs{B}=\kappa$. Fix a bijection $b:\kappa\to B$. For each $\xi<\kappa$, the set $b[\xi]$ is bounded in $\alpha$. Put
\[
\begin{aligned}
    \delta_\xi&=\sup_{\zeta<\xi}(b(\zeta)+1)<\alpha,\\
    h_\xi&=g_{\delta_\xi}\circ(b\restriction\xi):\xi\to2.
\end{aligned}
\]
Since $b$ is injective and $b[\xi]\subseteq\delta_\xi\le\delta_\eta$ for $\xi<\eta<\kappa$, the coherence of $\langle g_\delta:\delta<\alpha\rangle$ gives $h_\eta\restriction\xi=^*h_\xi$. By $\mathrm{CU}_2(\kappa)$, choose a map $h:\kappa\to2$ that uniformizes $\langle h_\xi:\xi<\kappa\rangle$. Define $g^B:B\to2$ by $g^B(b(\xi))=h(\xi)$ for $\xi<\kappa$, and put
\[
    g=g^B\cup
      \bigcup_{\xi<\kappa}
      \bigl(g_{c_\xi}\restriction(c_\xi\setminus B)\bigr).
\]
By the definition of $B$, $g$ is a well-defined map. Since $\langle c_\xi:\xi<\kappa\rangle$ is cofinal in $\alpha$, the domain of $g$ is $\alpha$. 

Fix $\eta<\kappa$ and put $F=\{\delta<c_\eta:g(\delta)\ne g_{c_\eta}(\delta)\}$. By the definition of $g$, we have $F\subseteq B$. For every $\xi<\kappa$, the relation $h\restriction\xi=^*h_\xi$ and the coherence of $\langle g_\delta:\delta<\alpha\rangle$ give
\[
    g\restriction(b[\xi]\cap c_\eta)
    =^*g_{\delta_\xi}\restriction(b[\xi]\cap c_\eta)
    =^*g_{c_\eta}\restriction(b[\xi]\cap c_\eta).
\]
Hence $b^{-1}(F)\cap\xi$ is finite for every $\xi<\kappa$. Since $\kappa$ is regular and uncountable, $b^{-1}(F)$ is finite, and hence so is $F$. Thus $g\restriction c_\eta=^*g_{c_\eta}$ for every $\eta<\kappa$, proving $\mathrm{CU}_2(\alpha)$.
\end{proof}

Combining Proposition~\ref{prop:cu-cofinality} with Theorem~\ref{thm:ordinal-rigidity}, we obtain the following corollary.

\begin{corollary}\label{cor:rigidity-cofinality}
Let $\alpha$ be a limit ordinal. Then $\alpha$ is rigid if and only if $\operatorname{cf}(\alpha)$ is rigid.
\end{corollary}

We recall some set-theoretic principles used in the next corollary. 

An ideal $\mathcal I$ on a set $X$ is called a \emph{$P$-ideal} if, for every sequence $\langle A_n:n<\omega\rangle$ in $\mathcal I$, there is $A\in\mathcal I$ such that $A_n\setminus A$ is finite for every $n<\omega$. The \emph{$P$-ideal dichotomy} ($\mathsf{PID}$) asserts that, for every set $X$ and every $P$-ideal $\mathcal I$ of countable subsets of $X$, one of the following holds:
\begin{enumerate}[label=\upshape(\roman*), leftmargin=*, widest=ii]
    \item There is an uncountable set $Y\subseteq X$ such that $[Y]^{\omega}\subseteq\mathcal I$.
    \item $X=\bigcup_{n<\omega}X_n$, where $A\cap X_n$ is finite for every $A\in\mathcal I$ and every $n<\omega$.
\end{enumerate}
$\mathsf{PID}$ is consistent with $\mathsf{ZFC}$ relative to the consistency of the existence of a supercompact cardinal \cite[Lemma~7]{todorcevic2000dichotomy}.

For a cardinal $\kappa$, the principle $\square(\kappa)$ asserts that there is a sequence $\langle C_\delta:\delta<\kappa\rangle$ satisfying the following conditions:
\begin{enumerate}[label=\upshape(\roman*), leftmargin=*, widest=iii]
    \item For every $\delta<\kappa$, the set $C_\delta$ is closed and unbounded in $\delta$.
    \item If $\beta$ is a limit point of $C_\delta$, then $C_\beta=C_\delta\cap\beta$.
    \item There is no closed unbounded set $C\subseteq\kappa$ such that $C_\beta=C\cap\beta$ for every limit point
    $\beta<\kappa$ of $C$.
\end{enumerate}
For an infinite cardinal $\lambda$, the principle $\square_\lambda$ asserts that there is a sequence $\langle C_\delta:\delta<\lambda^+\rangle$ satisfying conditions \textup{(i)} and \textup{(ii)} above with $\kappa=\lambda^+$, with the additional requirement that each $C_\delta$ have order type at most $\lambda$. In particular, $\square_\lambda$ implies $\square(\lambda^+)$
\cite[p.~110]{koenig2003local}.

\begin{corollary}\label{cor:small-cofinalities}
Let $\alpha$ be a limit ordinal.
\begin{enumerate}[label=\upshape(\roman*), leftmargin=*, widest=iii]
\item\label{it:cf-omega} If $\operatorname{cf}(\alpha)=\omega$, then $\alpha$ is rigid.
\item\label{it:cf-omega-one} If $\operatorname{cf}(\alpha)=\omega_1$, then $\alpha$ is not rigid.
\item\label{it:cf-pid} Assume $\mathsf{PID}$. If $\operatorname{cf}(\alpha)>\omega_1$, then $\alpha$ is rigid.
\item\label{it:cf-constructible} Assume $V=L$. If $\operatorname{cf}(\alpha)$ is an uncountable successor cardinal, then $\alpha$ is not rigid.
\end{enumerate}
\end{corollary}

\begin{proof}
Put $\kappa=\operatorname{cf}(\alpha)$. 

By Corollary~\ref{cor:rigidity-cofinality}, it suffices to determine whether $\kappa$ is rigid. The space $\omega$ is discrete and hence rigid. By Corollary~\ref{cor:omega-one-not-rigid}, $\omega_1$ is not rigid. This proves \ref{it:cf-omega} and \ref{it:cf-omega-one}.

For \ref{it:cf-pid}, assume $\mathsf{PID}$ and $\kappa>\omega_1$. By \cite[Theorem~3.16]{bergfalkLambieHanson} with $A=2$, $\mathrm{CU}_2(\kappa)$ holds. Thus $\kappa$ is rigid by Theorem~\ref{thm:ordinal-rigidity}.

For \ref{it:cf-constructible}, assume $V=L$ and let $\kappa$ be an uncountable successor cardinal. Since $L\models\mathsf{ZFC}$, the case $\kappa=\omega_1$ follows from \ref{it:cf-omega-one}. Otherwise, write $\kappa=\lambda^+$, where $\lambda$ is an uncountable cardinal. By  \cite[Theorem~27.1]{jech2003set}, $\square_\lambda$ holds, and hence so does $\square(\kappa)$. By \cite[Corollary~3.11]{koenig2003local}, there is a nontrivial coherent family $\langle f_\xi:\xi<\kappa\rangle$ with $f_\xi:\xi\to2$ for each $\xi<\kappa$. Thus $\mathrm{CU}_2(\kappa)$ fails, so $\kappa$ is not rigid by Theorem~\ref{thm:ordinal-rigidity}.
\end{proof}

\subsection*{Acknowledgements}We thank Hang Zhang of Southwest Jiaotong University for drawing our attention to Velichko's question on the rigidity of $\omega_1$, for helpful discussions, and for reading the manuscript carefully. We also thank an anonymous referee for a careful reading of an earlier version of this manuscript and for pointing out a gap concerning the use of the Axiom of Choice. The referee suggested that $\diamondsuit_{\omega_1}$ might imply that $\omega_1$ is not rigid. Pursuing this suggestion led us first to a proof under $\diamondsuit_{\omega_1}$ and later to a proof in $\mathsf{ZFC}$.

\subsection*{Declaration on the use of generative AI}
The authors used SnakeGPT (snake.gpt), a third-party API relay service, to access the ChatGPT-6.1 Sol model to improve the language and readability of the manuscript and to locate relevant references. The tool also helped the authors organize recurring arguments from an earlier draft into Lemmas~\ref{lem:section5-coding} and~\ref{lem:binary-realization}, thereby simplifying the proof of Theorem~\ref{thm:ordinal-rigidity}. The authors take full responsibility for the article.



\end{document}